\documentclass[submission]{FPSAC2019}

\usepackage{etex}
\reserveinserts{28}
\usepackage{pict2e}
\usepackage{xypic}
\usepackage{tikz}
\usepackage{MnSymbol}
\usepackage{latexsym}
\usepackage{stmaryrd}
\usepackage{pictex}
\usepackage{amsmath,amscd}
\usepackage{afterpage}
\usepackage{youngtab}
\usepackage{psfrag}
\usepackage[boxsize=.3 em]{ytableau}

\usepackage{verbatim}

\usepackage{graphics}
\usepackage{graphicx}

\font\co=lcircle10

\def\jr{\rotatedown{\smash{\raise2pt\hbox{\co \rlap{\rlap{\char'005} \char'007}}
               \raise6pt\hbox{\rlap{\vrule height6.5pt}}
                \raise2pt\hbox{\rlap{\hskip4pt \vrule
          height0.4pt depth0pt
                width7.7pt}}}}}
\def\textcross{\ \smash{\lower4pt\hbox{\rlap{\hskip4.15pt\vrule height14pt}}
                \raise2.8pt\hbox{\rlap{\hskip-3pt \vrule height.4pt depth0pt
                                width14.7pt}}}\hskip12.7pt}

\def\textelbow{\ \hskip.1pt\smash{\raise2.75pt%
                \hbox{\co \hskip 4.15pt\rlap{\rlap{\char'004} \char'006}
                \lower6.8pt\rlap{\vrule height3.5pt}
                \raise3.6pt\rlap{\vrule height3.5pt}}
                \raise2.8pt\hbox{%
                  \rlap{\hskip-7.15pt \vrule height.4pt depth0pt
width3.5pt}%
                  \rlap{\hskip4.05pt \vrule height.4pt depth0pt
width3.5pt}}}
                \hskip8.7pt}

\newtheorem{thm}{Theorem}[section]
\newtheorem{theorem}[thm]{Theorem}

\newtheorem{conjecture}[thm]{Conjecture}

\newtheorem{lemma}[thm]{Lemma}

\newtheorem{proposition}[thm]{Proposition}
\theoremstyle{definition}

\newtheorem{definition}[thm]{Definition}
\newtheorem{rem}[thm]{Remark}  
\newtheorem{remark}[thm]{Remark}  

\numberwithin{equation}{section}

\newcommand{\set}[1]{\left\{#1\right\}}

\newcommand{\CC}{\mathbb C}

\newcommand{\C}{\mathbb{C}}

\newcommand{\xx}{\widetilde{\mathbf{x}}}

\newcommand{\Gl}{\operatorname{GL}}
\newcommand{\Fl}{\operatorname{Fl}}

\def\Fcal{\mathcal{F}}

\def\TT{\mathbb{T}}
\def\Xcal{\mathcal{X}}
\def\Acal{\mathcal{A}}
\def\ZZ{\mathbb{Z}}
\newcommand{\overunder}[2]{
\!\begin{array}{c}
\scriptstyle{#1}\\[-.1in]
-\!\!\!-\!\!\!-\\[-.1in]
\scriptstyle{#2}
\end{array}
\!
}

\newcommand{\pathsw}[1]{L^\swarrow_{#1}}

\newcommand{\partne}[1]{\lambda^\nearrow(#1)} 
\newcommand{\partsw}[1]{\lambda^\swarrow(#1)}

\newcommand{\ppermsw}[1]{\pi^\swarrow_{#1}}

\title{Combinatorics of cluster structures in Schubert varieties}

\author{K. Serhiyenko\thanks{\href{mailto:khrystyna.serhiyenko@berkeley.edu}{khrystyna.serhiyenko@berkeley.edu} \addressmark{1}}, 
M. Sherman-Bennett\thanks{\href{mailto:m\_shermanbennett@berkeley.edu}{m\_shermanbennett@berkeley.edu} \addressmark{2}}, 
\and L. Williams\thanks{\href{mailto:williams@math.harvard.edu}{williams@math.harvard.edu} \addressmark{3}}}
\address{\addressmark{1}Department of Mathematics, University of California at Berkeley \\ \addressmark{2} Department of Mathematics, University of California at Berkeley \\ \addressmark{3} Department of Mathematics, Harvard University}

\received{\today}

\abstract{We give an explicit combinatorial description of 
cluster structures in Schubert varieties of the Grassmannian in terms of 
(target labelings of)
Postnikov's \emph{plabic graphs}.  This description is a natural generalization
of the description given by (Scott 2006) for the Grassmannian
	and has been
	believed by experts essentially since (Scott 2006), 
	though the statement was not formally written down
until (M\"uller-Speyer 2016).
To prove this conjecture we use a result of (Leclerc 2016), who used the module category of the preprojective algebra  
to prove that coordinate rings of 
many Richardson varieties in the complete flag variety 
admit cluster structures. We also adapt a construction of (Karpman 2016) 
to 
build cluster seeds associated to reduced expressions.
Further, we explicitly describe cluster structures 
in skew Schubert varieties using plabic graphs
whose boundary vertices need not be labeled in cyclic order.
\paragraph{R\'esum\'e.}
Nous donnons une description combinatoire explicite des structures 
en amas des vari\'et\'es de Schubert des  grassmanniennes en termes de graphes plabiques.
Cette description est une g\'en\'eralisation naturelle de la description donn\'ee par (Scott
2006) pour les grassmanniennes, et correspond \`a ce qui \'etait attendu par les experts depuis (Scott
2006), bien que l'\'enonc\'e  n'ait pas \'et\'e formellement \'ecrit avant (Muller-Speyer
2016).
Nous d\'ecrivons aussi explicitement les structures en amas des vari\'et\'es de Schubert g\'en\'eralis\'ees
en utilisant des graphes plabiques dont les sommets ne sont pas num\'erot\'es dans l'ordre cyclique.
}


\usepackage[backend=bibtex]{biblatex}
\begin{document}

\maketitle

\section{Introduction}

The main result of this extended abstract is that the coordinate ring of (the affine cone over)
any Schubert variety
of the Grassmannian (embedded into projective space via the Pl\"ucker embedding)
admits a cluster algebra structure, which is described explicitly in terms of 
plabic graphs. In the first section, we give some background on Schubert varieties, state our main result, and discuss some applications. Background material on cluster algebras and plabic graphs is delayed to \cref{sec:background}. 

 \emph{Cluster algebras} are a class of commutative rings
which were  introduced by Fomin and Zelevinsky \cite{ca1}; they are connected to 
many fields of mathematics including Teichm\"uller theory and quiver representations,
and their generators satisfy many nice properties, including the \emph{Laurent phenomenon}
\cite{ca1} 
and \emph{positivity theorem} \cite{LeeSchiffler, GHKK}. 
\emph{Plabic graphs} are certain planar bicolored
graphs which were introduced by Postnikov \cite{Postnikov}; plabic graphs (or rather an 
equivalent object, namely alternating strand diagrams) were subsequently used by 
Scott \cite{Scott} to show that the coordinate ring of the affine cone over the 
Grassmannian in its Pl\"ucker embedding admits a cluster algebra structure.

There is a natural plabic graph 
generalization of Scott's construction which experts have believed
for some time should give a cluster structure for Schubert varieties 
(and more generally,
positroid varieties).  
This construction was stated explicitly as a conjecture in 
a recent paper of M\"uller-Speyer \cite{MullerSpeyer0}, who 
additionally provided some 
evidence in \cite{MullerSpeyer}.  The conjecture 
can be stated roughly as follows.
\begin{conjecture}\label{conj:vague}
	Let $G$ be a reduced plabic graph corresponding to a Schubert (or more generally a positroid) variety.
	Then the target labeling of the faces of $G$ (which we identify
	with a collection of Pl\"ucker coordinates)
together with the dual quiver of $G$ gives rise to a seed for a cluster
	structure in the coordinate ring of the Schubert (or positroid) 
	variety.
\end{conjecture}
Meanwhile, Leclerc \cite{Leclerc} constructed a subcategory $\mathcal{C}_{v,w}$ of the module category of the preprojective algebra that has a cluster structure,
to show that the coordinate ring of 
each Richardson variety $\mathcal{R}_{v,w}$ of the 
complete flag variety contains a subalgebra which is a cluster algebra; when 
$w$ has a factorization of the form $w=xv$ with 
$\ell(w) = \ell(x)+\ell(v)$,  he showed that 
this subalgebra coincides with the 
coordinate ring.  
Because Schubert varieties are isomorphic to Richardson varieties with 
the above property,
Leclerc's result implies that their coordinate rings admit a cluster structure.
However, Leclerc's description of the cluster structure is 
very different from the 
plabic graph description and is far from explicit:
e.g. his cluster quiver is defined in terms of 
morphisms between modules of the preprojective algebra. 

In this paper we prove 
\cref{conj:vague} for Schubert varieties by relating 
Leclerc's cluster structure to the conjectural one coming from plabic graphs.
We also generalize our result to the setting of \emph{skew Schubert varieties}; 
 interestingly, these cluster structures for skew Schubert varieties
 depart from the one in 
\cref{conj:vague}, since they 
use \emph{generalized} plabic graphs (with 
boundary vertices which are no longer cyclically labeled).

Once we have proved that the coordinate rings of Schubert and skew Schubert
varieties have cluster structures, we obtain a number of results ``for free''
 from the cluster theory, including the Laurent phenomenon
and the positivity theorem for cluster variables.
We also obtain many combinatorially 
explicit cluster seeds for each Schubert and skew Schubert
variety. Other applications of our results, including a characterization of which Schubert varieties exhibit finite type cluster structures, are described in  \cref{sec:applications}.

\subsection{Notation for the flag variety}

Let $\Gl_n$ denote the general linear group, $B$ the Borel subgroup of 
lower triangular matrices, $B^+$ the opposite Borel subgroup of 
upper triangular matrices, and 
$W=S_n$ the Weyl group (which in this case is the symmetric group on $n$ letters).
$W$ is generated by the simple reflections $s_i$ for $1 \leq i \leq n-1$, where 
$s_i$ is the transposition exchanging $i$ and $i+1$, and it contains 
a longest element, which we denote by $w_0$, with $\ell(w_0) = {n \choose 2}$.
The \emph{complete flag variety} $\Fl_n$ is the homogeneous
space $B\setminus \Gl_n$. 
Concretely, each element $g$ of $\Gl_n$ gives rise to a flag of 
subspaces $\{V_1 \subset V_2 \subset \dots \subset V_n\}$, 
where $V_i$ denotes the span of the top $i$ rows of $g$.
The action of $B$ on the left preserves the flag, so we can identify 
$\Fl_n$ with 
the set of \emph{flags} $\{V_1 \subset V_2 \subset \dots \subset V_n\}$
where $\dim V_i = i$.

Let $\pi:\Gl_n \to \Fl_n$  
denote the natural projection 
$\pi(g):= Bg$.  The Bruhat decomposition 
$$\Gl_n = \bigsqcup_{w\in W} BwB$$
projects to the Schubert decomposition 
$$\Fl_n = \bigsqcup_{w\in W} C_w$$
where $C_w = \pi(BwB)$ is the \emph{Schubert cell} associated to $w$, 
 isomorphic to $\C^{\ell(w)}$.
We also have the Birkhoff decomposition 
$$\Gl_n = \bigsqcup_{w\in W} BwB^+,$$ which projects to the 
{opposite} Schubert decomposition 
$$\Fl_n = \bigsqcup_{w\in W} C^w$$
where $C^w = \pi(B w B^+)$ is the \emph{opposite Schubert cell} 
associated to $w$, isomorphic to $\C^{\ell(w_0) - \ell(w)}$.

The intersection $$\mathcal{R}_{v,w}:=C^v \cap C_w$$
has been considered by Kazhdan and Lusztig \cite{KL} in relation to 
Kazhdan-Lusztig polynomials.  $\mathcal{R}_{v,w}$ is 
nonempty only if $v \leq w$ in the Bruhat order of $W$, 
and it is then a smooth irreducible locally closed subset of $C_w$ 
of dimension $\ell(w)-\ell(v)$.  Sometimes $\mathcal{R}_{v,w}$
is called an \emph{open Richardson variety} \cite{KLS} because its closure 
is a \emph{Richardson variety} \cite{Rich}.
We have a stratification of the complete flag variety
$$\Fl_n = \bigsqcup_{v \leq w} \mathcal{R}_{v,w}.$$



\subsection{Notation for 
 the Grassmannian}

Fix $1 < k < n$.
The parabolic subgroup 
$W_K = \langle s_1,\dots,s_{k-1}\rangle \times \langle s_{k+1}, s_{k+2},\dots,
s_{n-1} \rangle < W$ 
gives rise to a parabolic subgroup 
$P_K = \bigsqcup_{w \in W_K} B \dot{w} B$ in $\Gl_n$, where $\dot{w}$ is a matrix representative for $w$
in $\Gl_n$.
The longest element $w_K$ of $W_K$ has length
 $\ell(w_K) = {k \choose 2} + {n-k \choose 2}$. 

The \emph{Grassmannian} $Gr_{k,n}$ 
is the homogeneous space 
$P_K\setminus \Gl_n$.
We can think of the 
{Grassmannian} $Gr_{k,n} = 
P_K\setminus \Gl_n$ more concretely as 
the set of all $k$-planes in an $n$-dimensional
vector space $\CC^n$. 
An element of
$Gr_{k,n}$ can be viewed as a full rank 
$k\times n$ matrix of rank $k$, modulo left
multiplication by invertible $k\times k$ matrices.  
That is, two $k\times n$ matrices of rank $k$ represent the same point in $Gr_{k,n}$ if and only if they
can be obtained from each other by invertible row operations.

For a positive integer $a$ and $w \in W$, let $[a]:=[1, a]$ and $w[a]:=\{w(1), \dots, w(a)\}$. Let $\binom{[n]}{k}$ be the set of all $k$-element 
subsets of $[n]$. 

Given $V\in Gr_{k,n}$ represented by a $k\times n$ matrix $A$, for $I\in \binom{[n]}{k}$ we let $\Delta_I(V)$ be the maximal minor of $A$ located in the column set $I$. The $\Delta_I(V)$ do not depend on our choice of matrix $A$ (up to simultaneous rescaling by a nonzero constant), and are called the {\itshape Pl\"{u}cker coordinates} of $V$.
The Pl\"ucker coordinates give an embedding of $Gr_{k,n}$ into projective space of dimension
${n \choose k} - 1$.

We have the usual projection $\pi_k$ from the complete flag variety $\Fl_n$ to the 
Grassmannian $Gr_{k,n}$.
Let $W^K = W^K_{\min}$ and 
$W^K_{\max}$ denote the set of minimal- and
maximal-length coset representatives for $W_K \setminus W$. 

Rietsch studied the projections 
of the open Richardson varieties in the complete flag variety to 
partial flag varieties \cite{RietschThesis}.  In particular,
when  $v\in W^K_{\max}$
 (or when $w\in W^K_{\min}$), 
the projection $\pi_k$ is an isomorphism from $\mathcal{R}_{v,w}$ to 
 $\pi_k(\mathcal{R}_{v,w})$.  
We obtain a stratification 
$$Gr_{k,n}=\bigsqcup_{v\leq w} \pi_k(\mathcal{R}_{v,w})$$
where $(v,w)$ range over all $v\in W^K_{\max}$, $w\in W$, such that 
$v\leq w$.
 Following work of Postnikov \cite{Postnikov, KLS},
 the strata $\pi_k(\mathcal{R}_{v,w})$ are sometimes called 
 \emph{open positroid varieties}, while their closures
 are called \emph{positroid varieties}.

It follows from the definitions (see e.g. \cite[Section 6]{KLS}) that 
positroid varieties include Schubert varieties
in the Grassmannians, which we now define.

\begin{definition}
Let $I\in \binom{[n]}{k}$.
The \emph{Schubert cell} $\Omega_{I}$  is defined to be the set
$$ \{A \in Gr_{k,n}\ \vert \ \text{the lexicographically minimal nonvanishing
Pl\"ucker coordinate of }A\text{ is }\Delta_{I}(A) \}.$$
	The \emph{Schubert variety} $X_{I}$ is the closure
	$\overline{\Omega_I}$ of $\Omega_I$.

\end{definition}

When $v\in W^K_{\max}$,
$\overline{\pi_k(\mathcal{R}_{v,w_0})}$ is isomorphic to the 
Schubert variety $X_{v^{-1}[k]}$ in the Grassmannian, 
which has dimension 
$\ell(w_0)-\ell(v)$. 
More generally, if $v\in W^K_{max}$ and $w\in W$ 
has a factorization of the form $w=xv$ which is \emph{length-additive}, i.e. where 
$\ell(w) = \ell(x)+\ell(v)$, then we refer to 
$\overline{\pi_k(\mathcal{R}_{v,w})}$ as a \emph{skew Schubert variety}.

Let $\lambda$
 denote a Young diagram contained in a $k \times (n-k)$ rectangle.  
We can identify $\lambda$ with its southeast boundary; we think of the boundary as a lattice path $\pathsw{\lambda}$ in the rectangle from the northeast corner to the southwest corner taking steps west 
and south
.
Labeling the steps of $\pathsw{\lambda}$ from $1$ to $n$, the labels of the 
south steps give a $k$-element subset of $[n]$.
Conversely, each $I \in \binom{[n]}{k}$ can be identified with a 
Young diagram, denoted $\partsw{I}$. 
(Later we will also need an analogous partition $\partne{I}$ obtained by reading the boundary of $\lambda$ from 
southwest to northeast.)
The map $\partsw$  allows us to index Schubert 
cells and varieties by Young diagrams, denoting them 
$\Omega_{\lambda}$ and
$X_{\lambda}$, 
respectively.  
The dimension of $\Omega_{\lambda}$ and $X_{\lambda}$ is $|\lambda|$, the number of boxes of 
$\lambda$.



\subsection{The main result}

We now state the main result.  Note that 
the definitions of plabic graph and trip permutation can be found in 
\cref{sec:background}.

We associate with a Young diagram $\lambda$ a permutation
$\ppermsw{\lambda}$: 
in list notation, this permutation is obtained by first reading 
the labels of the horizontal steps of $\pathsw{\lambda}$, and then reading the labels of the 
vertical steps of $\pathsw{\lambda}$.
(Moreover any fixed points in positions $1,2,\dots, n-k$ are ``black"
and any fixed points in positions $n-k+1,\dots, n$ are ``white.")

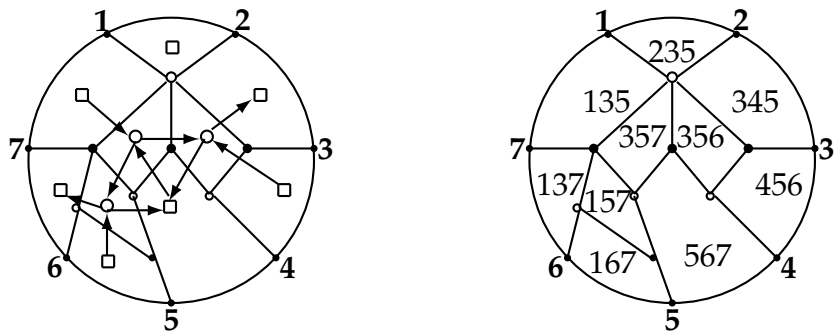
\begin{figure}[h]
\begin{center}
\setlength{\unitlength}{.9pt}
\begin{picture}(100,120)(-70,-70)
\thicklines
\put(15,30){\circle{4}}
\put(-18,0){\circle*{4}}
\put(15,0){\circle*{4}}
\put(47,0){\circle*{4}}
\put(13,30){\line(-31,23){25}}
\put(17,30){\line(31,23){25}}
\put(13,28){\line(-25,-23){30}}
\put(17,28){\line(25,-23){30}}
\put(15,28){\line(0,-1){29}}
\put(-18,0){\line(-1,0){28}}
\put(-18,0){\line(1,-1.1){17}}
\put(47,0){\line(-1,-1.2){16}}
\put(47,0){\line(1,0){28}}
\put(15,-5){\circle{120}}
\put(15,-65){\circle*{3}}
\put(-1,-20){\circle{3}}
\put(-12,48){\circle*{3}}
\put(31,-20){\circle{3}}
\put(42,48){\circle*{3}}
\put(-45,0){\circle*{3}}
\put(75,0){\circle*{3}}
\put(-25,-25){\circle{3}}
	\put(7,-46){\circle*{3}}
\put(-29,-46){\circle*{3}}
\put(59,-46){\circle*{3}}
\put(-29,-46){\line(5,20){5}}
	\put(-24,-25){\line(5,20){6.5}}
	\put(-24,-25){\line(10,-7){30}}
	\put(-1,-20){\line(5,-14){16}}
	\put(15,0){\line(-10,-12){16}}
	\put(15,0){\line(10,-12){16}}
\put(59,-46){\line(-20,20){26}}         
\put(-14,54){\makebox(0,0){$\mathbf{1}$}}
\put(45,54){\makebox(0,0){$\mathbf{2}$}}
\put(80,0){\makebox(0,0){$\mathbf{3}$}}
\put(64,-50){\makebox(0,0){$\mathbf{4}$}}
\put(15,-72){\makebox(0,0){$\mathbf{5}$}}
\put(-34,-50){\makebox(0,0){$\mathbf{6}$}}
\put(-50,0){\makebox(0,0){$\mathbf{7}$}}

\put(13,40){\line(0,1){5}}
\put(13,40){\line(1,0){5}}
\put(13,45){\line(1,0){5}}
\put(18,40){\line(0,1){5}}

\put(50,20){\line(0,1){5}}
\put(50,20){\line(1,0){5}}
\put(50,25){\line(1,0){5}}
\put(55,20){\line(0,1){5}}

\put(60,-20){\line(0,1){5}}
\put(60,-20){\line(1,0){5}}
\put(60,-15){\line(1,0){5}}
\put(65,-20){\line(0,1){5}}

\put(-34,-20){\line(0,1){5}}
\put(-34,-20){\line(1,0){5}}
\put(-34,-15){\line(1,0){5}}
\put(-29,-20){\line(0,1){5}}

\put(-14,-50){\line(0,1){5}}
\put(-14,-50){\line(1,0){5}}
\put(-14,-45){\line(1,0){5}}
\put(-9,-50){\line(0,1){5}}

\put(12,-27){\line(0,1){5}}
\put(12,-27){\line(1,0){5}}
\put(12,-22){\line(1,0){5}}
\put(17,-27){\line(0,1){5}}

\put(-25,20){\line(0,1){5}}
\put(-25,20){\line(1,0){5}}
\put(-25,25){\line(1,0){5}}
\put(-20,20){\line(0,1){5}}

\put(30,5){\circle{5}}
\put(0,5){\circle{5}}
\put(-12,-24){\circle{5}}
\put(32,8){{\vector(4,3){17}}}
\put(29,4){{\vector(-4,-7){15}}}
\put(2,4){{\vector(1,0){25}}}
\put(0,3){{\vector(-5,-10){12}}}
\put(-12,-26){{\vector(1,0){24}}}
\put(-14,-25){{\vector(-3,1){16}}}
\put(-12,-45){{\vector(0,1){18}}}
\put(14,-20){{\vector(-2,3){15}}}
\put(-20,20){{\vector(15,-13){17}}}
\put(60,-15){{\vector(-3,2){28}}}
\end{picture}
\qquad \qquad
\qquad \qquad
\begin{picture}(100,120)(-70,-70)
\thicklines
\put(15,30){\circle{4}}
\put(-18,0){\circle*{4}}
\put(15,0){\circle*{4}}
\put(47,0){\circle*{4}}
\put(13,30){\line(-31,23){25}}
\put(17,30){\line(31,23){25}}
\put(13,28){\line(-25,-23){30}}
\put(17,28){\line(25,-23){30}}
\put(15,28){\line(0,-1){29}}
\put(-18,0){\line(-1,0){28}}
\put(-18,0){\line(1,-1.1){17}}
\put(47,0){\line(-1,-1.2){16}}
\put(47,0){\line(1,0){28}}
\put(15,-5){\circle{120}}
\put(15,-65){\circle*{3}}
\put(-1,-20){\circle{3}}
\put(-12,48){\circle*{3}}
\put(31,-20){\circle{3}}
\put(42,48){\circle*{3}}
\put(-45,0){\circle*{3}}
\put(75,0){\circle*{3}}
\put(-25,-25){\circle{3}}
	\put(7,-46){\circle*{3}}
\put(-29,-46){\circle*{3}}
\put(59,-46){\circle*{3}}
\put(-29,-46){\line(5,20){5}}
	\put(-24,-25){\line(5,20){6.5}}
	\put(-24,-25){\line(10,-7){30}}
	\put(-1,-20){\line(5,-14){16}}
	\put(15,0){\line(-10,-12){16}}
	\put(15,0){\line(10,-12){16}}
\put(59,-46){\line(-20,20){26}}         
\put(-14,54){\makebox(0,0){$\mathbf{1}$}}
\put(44,54){\makebox(0,0){$\mathbf{2}$}}
\put(80,0){\makebox(0,0){$\mathbf{3}$}}
\put(64,-50){\makebox(0,0){$\mathbf{4}$}}
\put(15,-72){\makebox(0,0){$\mathbf{5}$}}
\put(-34,-50){\makebox(0,0){$\mathbf{6}$}}
\put(-50,0){\makebox(0,0){$\mathbf{7}$}}

\put(15,40){\makebox(0,2){$235$}}
\put(15,15){\makebox(24,-20){$356$}}
\put(15,15){\makebox(-25,-20){$357$}}
\put(0,30){\makebox(-25,-20){$135$}}
\put(50,20){\makebox(0,0){$345$}}
\put(60,-15){\makebox(0,0){$456$}}
\put(30,-45){\makebox(0,0){$567$}}
\put(-10,-48){\makebox(0,0){$167$}}
\put(-32,-15){\makebox(0,0){$137$}}
\put(-12,-22){\makebox(0,0){$157$}}
\end{picture}
\end{center}
	\caption{A plabic graph $G$ for $Gr_{3,7}$ with trip permutation 
	$\ppermsw{\lambda}=(2,4,6,7,1,3,5)$, for $\lambda=(4, 3, 2)$, together with the dual quiver of $G$ and the target face 
	labeling.}
\label{fig:plabic2}
\end{figure}

\begin{theorem}\label{thm:main}
	Consider the Schubert variety $X_{\lambda}$ of $Gr_{k,n}$. 
	Let $G$ be a reduced plabic graph (with boundary vertices
	labeled clockwise from $1$ to $n$) with trip permutation 
	$\ppermsw{\lambda}$.
Construct the dual quiver of $G$ and label its vertices by the Pl\"ucker coordinates
given by the target labeling of $G$, see \cref{def:faces} and \cref{fig:plabic2}.
The coordinate ring of (the affine cone over) $X_{\lambda}$ is a 
cluster algebra, and this labeled quiver gives a seed for this cluster algebra.
\end{theorem}


We actually prove something a bit more general than \cref{thm:main}.

\begin{theorem}\label{thm:main2}
Consider the skew Schubert variety 
	$\overline{\pi_k(\mathcal{R}_{v,w})}$,
where $v\in W^K_{max}$ and 
$w$ has a length-additive factorization $w=xv$.
Let $G$ be a reduced plabic graph
(with boundary vertices labeled clockwise from $1$ to $n$) with 
trip permutation $vw^{-1} = x^{-1}$,  such that boundary lollipops
are white if and only if they are in $[k]$. Then the coordinate
ring of (the affine cone over) the 
skew Schubert variety $\overline{\pi_k(\mathcal{R}_{v,w})}$
is a cluster algebra, and $G$ gives rise to a seed as follows: apply $v^{-1}$ to the boundary vertices of $G$ and then label the dual quiver using the target labeling.
\end{theorem}

In the case of Schubert varieties, \cref{thm:main} resolves
\cref{conj:vague}.
Note that 
there is another version of the conjecture which uses the \emph{source
labeling} of $G$ instead of the 
target labeling \cite[Remark 3.5]{MullerSpeyer0}.
Both conjectures make sense more generally for positroid varieties
and arbitrary reduced plabic graphs (whose trip permutations can be
arbitrary decorated permutations).
However, the cluster structure that we give in \cref{thm:main2}
is different from either of the cluster structures proposed in 
\cite{MullerSpeyer0}.  

\begin{rem} When $\overline{\pi_k(\mathcal{R}_{v, w})}$ is not a skew Schubert variety, the seeds in Leclerc's cluster subalgebra 
in general do not come from generalized plabic graphs.  Indeed, for $v=(2, 5, 1, 4, 3)$ and $w=(5, 3, 4, 2, 1)$, the unique
	seed in Leclerc's cluster subalgebra for $\pi_k(\mathcal{R}_{v, w})$ consists of frozen variables and has extended cluster $\{\Delta_{13}, \Delta_{23}, \Delta_{14}, \Delta_{15}, \Delta_{45}\}$. These Pl\"ucker coordinates cannot be obtained as the set of face labels of any (generalized) plabic graph
	(note that the label $2$ occurs only once in the set of Pl\"ucker coordinates).
\end{rem}

\subsection{Outline of the proof}

While each skew Schubert variety 
$\overline{\pi_k(\mathcal{R}_{v,w})}$
(where 
$v=w_K v'\in W^K_{max}$ 
and $w\in W$ has a length-additive factorization 
$\mathbf{w}=\mathbf{xv}=\mathbf{x w_K v'}$ into reduced expressions for $x$, $w_K$, and $v'$)
corresponds to an equivalence class of plabic graphs, 
there is 
one among them which is particularly nice, which we call the \emph{rectangles seed}. 
The first step of our proof is to give an explicit description of the rectangles seed for each skew Schubert variety
$\overline{\pi_k(\mathcal{R}_{v,w})}$, in terms of $v$ and a Young diagram associated to $(v,w)$.

A construction of Karpman \cite{Karpman} produces
a bridge-decomposable plabic graph associated to a pair $(y,\mathbf{z})$, where 
$y^{-1} \in W^K_{\max}$, $\mathbf{z}$ is a reduced decomposition for $z$, and $y\leq z$.
The second step of our proof is to show that if we perform her construction for the pair 
$(w_K, \mathbf{x w_K})$ and then relabel boundary vertices of the resulting plabic graph by $v^{-1}$ to obtain a graph $G$, 
the target labeling of the dual quiver of $G$ gives rise to 
 the rectangles seed 
for $\overline{\pi_k(\mathcal{R}_{v,w})}$. 

 In \cite{Leclerc}, Leclerc produces a cluster seed
 associated to each pair $(v,\mathbf{w})$, where $v\in W^K_{max}$ and $v \leq w$.
 The third step of our proof is to verify that for the choice $(v, \mathbf{w} = \mathbf{x w_K v'})$,  
Leclerc's construction gives rise to the rectangles seed.
To prove \cref{thm:main2} from the previous steps, we show that mutations of the plabic graph $G$ (known as ``square moves") coincide with certain mutations of the rectangles seed. \cref{thm:main} is then deduced from \cref{thm:main2}.


\subsection{Applications of the main result}
\label{sec:applications}

In this section we sketch some applications of the main result,
including ring-theoretic properties for the coordinate rings of 
skew Schubert varieties, canonical bases, and a classification of the 
finite type cluster structures we obtain.

Combining \cref{thm:main} and \cref{thm:main2} with results of Muller and Speyer
\cite{Muller}, \cite[Theorem 3.3]{MullerSpeyer0}, 
we find that the coordinate rings of Schubert and skew Schubert varieties (viewed as cluster algebras) are \emph{locally acyclic},
which implies that they are finitely generated, normal, locally a complete
intersection, and equal to their own upper cluster algebras.

Combining our results with results of Ford-Serhiyenko \cite[Theorem 1.2]{FordSer}, 
we find that the quivers giving
rise to the cluster structures for Schubert and skew Schubert varieties admit
\emph{green-to-red sequences}, which by Gross-Hacking-Keel-Kontsevich \cite{GHKK} implies that the cluster algebras have
\emph{Enough Global Monomials} and hence 
each coordinate ring has a canonical basis of
\emph{theta functions}, parameterized by the lattice of $g$-vectors.


%
%

In \cite{Scott}, Scott classified the Grassmannians whose coordinate rings
have a cluster algebra of finite type.  He showed that in general the
cluster algebras have infinite type, except in the following cases:
the coordinate ring of $Gr_{2,n}$ is a cluster algebra of type $A_{n-3}$,
while the coordinate rings of $Gr_{3,6}$, $Gr_{3,7},$ and $Gr_{3,8}$ are
cluster algebras of types $D_4$, $E_6$, and $E_8$, respectively.

It is straightforward to classify for which skew Schubert varieties $\overline{\pi_k(\mathcal{R}_{v, w})}$ 
the cluster structure described here is finite type. It depends only on $wv^{-1}$.

\begin{proposition} Let $v \leq w$, where $v \in W^K_{\max}$ and $w=xv$ is length-additive. Let $\lambda=\partne{x[k]}$ 
	and let $\lambda'$ be the Young diagram obtained from $\lambda$ by removing all boxes that touch the southeast boundary of $\lambda$ either along an edge or at the southeast corner. Then the cluster structure on the coordinate ring of 
	$\overline{\pi_k(\mathcal{R}_{v, w}})$ given in \cref{thm:main2} is

\begin{enumerate}
\item type $A$ if and only if $\lambda'$ does not contain a $2 \times 2$ rectangle;
\item type $D$ if and only if $\lambda'=(i, 2)$ or its transpose for $i \geq 2$;
\item type $E_6$, $E_7$, or $E_8$ if and only if $\lambda'$  or its transpose is one of $(3, 3)$, $(3, 2, 1)$, $(4, 3)$, $(4, 2, 1)$, $(3, 3, 1)$, $(5, 3)$, $(5, 2, 1)$, $(4, 4)$, $(4, 2, 2)$.
\end{enumerate}
 
In particular, the Schubert variety $X_{\lambda}$ is of finite type if and only if $\lambda'$ is in the above list.
\end{proposition}




\section{Background on cluster structures and plabic graphs}
\label{sec:background}
\subsection{Background on cluster structures}
Cluster algebras are a class of rings with a particular 
combinatorial structure; they were introduced by Fomin and Zelevinsky in \cite{ca1}.

\begin{definition}
\label{quiver}
A \emph{quiver} $Q$ is a directed graph; we will assume that $Q$ has no 
loops or $2$-cycles.
Each vertex is designated either  \emph{mutable} or \emph{frozen}.
\end{definition}

\begin{definition}
Let $q$ be a mutable vertex of quiver $Q$.  The quiver mutation 
$\mu_q$ transforms $Q$ into a new quiver $Q' = \mu_q(Q)$ via a sequence of three steps:
\begin{enumerate}
\item For each oriented two path $r \to q \to s$, add a new arrow $r \to s$
(unless $r$ and $s$ are both frozen, in which case do nothing).
\item Reverse the direction of all arrows incident to the vertex $q$.
\item Repeatedly remove oriented $2$-cycles until unable to do so.
\end{enumerate}
\end{definition}

\begin{definition}
\label{def:seed0}
Choose $M\geq N$ positive integers.
Let $\Fcal$ be an \emph{ambient field}
of rational functions
in $N$ independent
variables
over
$\C(x_{N+1},\dots,x_M)$.
A \emph{labeled seed} in~$\Fcal$ is
a pair $(\xx, Q)$, where
$\xx = (x_1, \dots, x_M)$ forms a free generating
set for
$\Fcal$,
and
$Q$ is a quiver on vertices
$1, 2, \dots,N, N+1, \dots, M$,
whose vertices $1,2, \dots, N$ are
\emph{mutable}, and whose vertices $N+1,\dots, M$ are \emph{frozen}.

We refer to~$\xx$ as the (labeled)
\emph{extended cluster} of a labeled seed $(\xx, Q)$.
The variables $\{x_1,\dots,x_N\}$ are called \emph{cluster
variables}, and the variables $c=\{x_{N+1},\dots,x_M\}$ are called
\emph{frozen} or \emph{coefficient variables}.
We often view the labeled seed as a quiver $Q$ where each vertex $i$ is labeled by 
the corresponding variable $x_i$.
\end{definition}

\begin{definition}
\label{def:seed-mutation0}
Let $(\xx, Q)$ be a labeled seed in $\Fcal$,
and let $q \in \{1,\dots,N\}$.
The \emph{seed mutation} $\mu_q$ in direction~$q$ transforms
$(\xx, Q)$ into the labeled seed
$\mu_q(\xx,  Q)=(\xx', \mu_q(Q))$, where the cluster
$\xx'=(x'_1,\dots,x'_M)$ is defined as follows:
$x_j'=x_j$ for $j\neq q$,
whereas $x'_q \in \Fcal$ is determined
by the \emph{exchange relation}
\begin{equation*}
\label{exchange relation0}
x'_q\ x_q = 
\prod_{q \to r} x_r + \prod_{s \to q} x_s,
\end{equation*}
where the first product is over all arrows $q \to r$ in $Q$
which start at $q$, and the second product is over all arrows $s\to q$ which 
end at $q$.
\end{definition}

\begin{remark}
It is not hard to check that seed mutation is an involution.
\end{remark}

\begin{remark}
Note that arrows between two frozen vertices of a quiver do not
affect seed mutation (they do not affect the mutated quiver
or the exchange relation).  For that reason, one may omit
arrows between two frozen vertices. 
\end{remark}

\begin{definition}
\label{def:patterns0}
Consider the \emph{$N$-regular tree}~$\TT_N$
whose edges are labeled by the numbers $1, \dots, N$,
so that the $N$ edges emanating from each vertex receive
different labels.
A \emph{cluster pattern}  is an assignment
of a labeled seed $\Sigma_t=(\xx_t, Q_t)$
to every vertex $t \in \TT_N$, such that the seeds assigned to the
endpoints of any edge $t \overunder{q}{} t'$ are obtained from each
other by the seed mutation in direction~$q$.
The components of
$\xx_t$ are written as $\xx_t = (x_{1;t}\,,\dots,x_{N;t}).$
\end{definition}

Clearly, a cluster pattern  is uniquely determined
by an arbitrary  seed.

\begin{definition}
\label{def:cluster-algebra0}
Given a cluster pattern, we denote
\begin{equation*}
\label{eq:cluster-variables0}
\Xcal
= \bigcup_{t \in \TT_N} \xx_t
= \{ x_{i,t}\,:\, t \in \TT_N\,,\ 1\leq i\leq N \} \ ,
\end{equation*}
the union of clusters of all the seeds in the pattern.
The elements $x_{i,t}\in \Xcal$ are called \emph{cluster variables}.
The
\emph{cluster algebra} $\Acal$ associated with a
given pattern is the $\ZZ[c]$-subalgebra of the ambient field $\Fcal$
generated by all cluster variables: $\Acal = \ZZ[c] [\Xcal]$.
We denote $\Acal = \Acal(\xx,  Q)$, where
$(\xx,Q)$
is any seed in the underlying cluster pattern.
In this generality,
$\Acal$ is called a \emph{cluster algebra from a quiver}, or a
\emph{skew-symmetric cluster algebra of geometric type.}
We say that $\Acal$ has \emph{rank $N$} because each cluster contains
$N$ cluster variables.
\end{definition}

\subsection{Background on plabic graphs}

In this section 
we review Postnikov's notion of \emph{plabic graphs} \cite{Postnikov}.

\begin{definition}
A {\it plabic (or planar bicolored) graph\/}
is an undirected graph $G$ drawn inside a disk
(considered modulo homotopy)
with $n$ {\it boundary vertices\/} on the boundary of the disk,
labeled $1,\dots,n$ in clockwise order, as well as some
colored {\it internal vertices\/}.
These internal vertices
are strictly inside the disk and are
colored in black and white. 
An internal vertex of degree one adjacent to a boundary vertex is a \emph{lollipop}.
We will always assume that no vertices of the same color are adjacent, and that 
each boundary vertex $i$ is adjacent to a single internal vertex.

\end{definition}

See \cref{fig:plabic2} for an example of a plabic graph.

There is a natural set of local transformations (moves) of plabic graphs, which we now describe.
Note that we will always assume that a plabic graph $G$ has no isolated 
components (i.e. every connected component contains at least
one boundary vertex).  We will also assume that $G$ is \emph{leafless}, 
i.e.\ if $G$ has an 
internal vertex of degree $1$, then that vertex must be adjacent to a boundary
vertex.

(M1) SQUARE MOVE (Urban renewal).  If a plabic graph has a square formed by
four trivalent vertices whose colors alternate,
then we can switch the
colors of these four vertices (and add some degree $2$ vertices to preserve
the bipartiteness of the graph).

(M2) CONTRACTING/EXPANDING A VERTEX.
Any degree $2$ internal vertex not adjacent to the boundary can be deleted,
and the two adjacent vertices merged.
This operation can also be reversed.  Note that this operation can always be used
to change an arbitrary
 square face of $G$ into a square face whose four vertices are all trivalent.

(M3) MIDDLE VERTEX INSERTION/REMOVAL.
We can always remove or add degree $2$ vertices at will, subject to the 
  condition that the graph remains bipartite.

\begin{figure}[h]
\centering
\includegraphics[height=.5in]{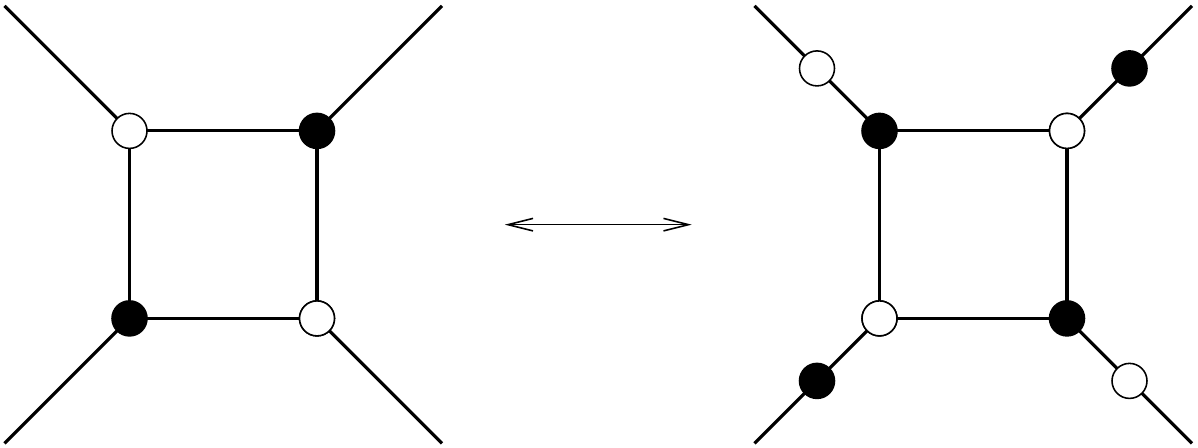}
\hspace{.5in}
\raisebox{6pt}{\includegraphics[height=.4in]{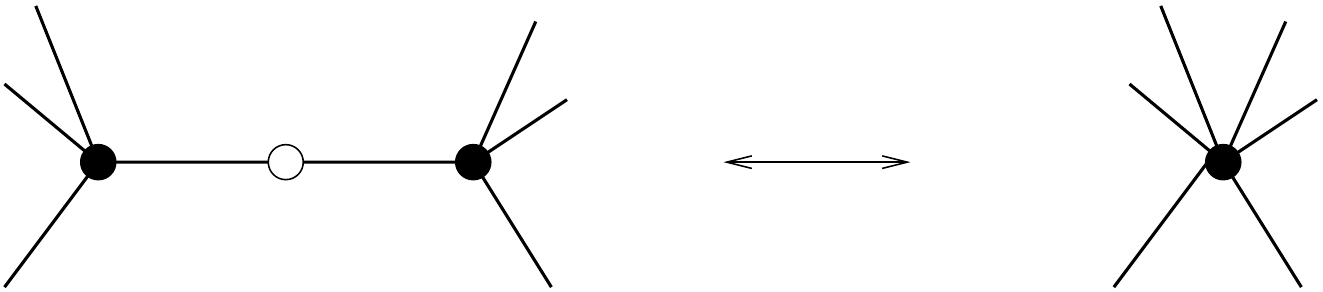}}
\hspace{.5in}
\raisebox{16pt}{\includegraphics[height=.07in]{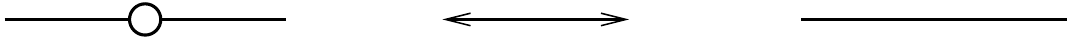}}
\caption{%
	Moves (M1), (M2), (M3).}
\label{M1}
\end{figure}

See \cref{M1} for depictions of these three moves.

(R1) PARALLEL EDGE REDUCTION.  If a plabic graph contains
two trivalent vertices of different colors connected
by a pair of parallel edges, then we can remove these
vertices and edges, and glue the remaining pair of edges together.

\begin{figure}[h]
\centering
\includegraphics[height=.25in]{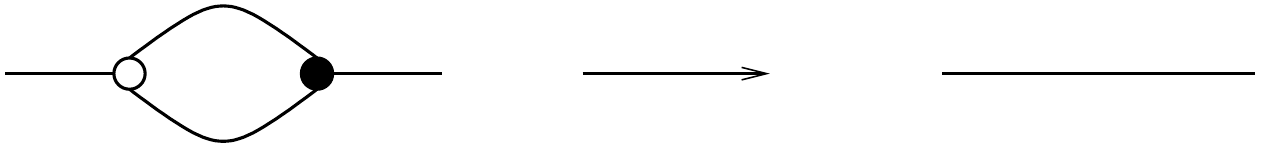}
\caption{Parallel edge reduction}
\label{R1}
\end{figure}

\begin{definition}
Two plabic graphs are called \emph{move-equivalent} if they can be obtained
from each other by moves (M1)-(M3).  The \emph{move-equivalence class}
of a given plabic graph $G$ is the set of all plabic graphs which are move-equivalent
to $G$.
A leafless plabic graph without isolated components
is called \emph{reduced} if there is no graph in its move-equivalence
class to which we can apply (R1).
\end{definition}

\begin{definition} A \emph{decorated permutation} $\pi^:$ is a permutation $\pi \in S_n$ together with a coloring $\set{i \ \vert \ \pi(i)=i} \to \{\text{black, white}\}$. 
\end{definition}

\begin{definition}\label{def:rules}
Given a reduced plabic graph $G$,
a  \emph{trip} $T$ is a directed path which starts at
some boundary vertex
$i$, and follows the ``rules of the road": it turns (maximally) right at a
black vertex,  and (maximally) left at a white vertex.
 Note that $T$ will also
end at a boundary vertex $j$; we then refer to this trip as
$T_{i \to j}$. Setting $\pi(i)=j$ for each such trip,
 we associate a (decorated) \emph{trip permutation}
$\pi_G=(\pi(1),\dots,\pi(n))$ to each reduced plabic graph $G$, where a fixed point $\pi(i)=i$ is colored white (black) if there is a white (black) lollipop at boundary vertex $i$.
\end{definition}

As an example, the trip permutation associated to the
reduced plabic graph in \cref{fig:plabic2} is $(2, 4, 6, 7, 1, 3, 5)$.

\begin{remark}\label{rem:moves}
Note that the trip permutation of a plabic graph is preserved 
by the local moves (M1)-(M3), but not by (R1). For reduced plabic graphs the converse holds, namely 
it follows from \cite[Theorem 13.4]{Postnikov} 
that any two reduced plabic graphs with the same trip permutation are 
move-equivalent.
\end{remark}

Now we use the notion of trips to label each face of $G$
by a Pl\"ucker coordinate.
Towards this end, note that every trip
will partition the faces of a plabic graph into
two parts: those on the left of the trip, and those on the right
of a trip.

\begin{definition}\label{def:faces}
Let $G$ be a reduced plabic graph with $b$ boundary vertices.
For each one-way trip $T_{i\to j}$ with $i \neq j$, we place the label $j$
 in every face which is to the left of $T_{i\to j}$. If $i=j$ (that is, $i$ is adjacent to a lollipop), we place the label $i$ 
in all faces if the lollipop is white and in no faces if the lollipop is black.
We then obtain a labeling $\mathcal{F}(G)$
of faces of $G$ by subsets of $[b]$ which
we call the \emph{target}
\emph{labeling} of $G$.  We identify each $a$-element subset of $[b]$
with the corresponding Pl\"ucker coordinate.
\end{definition}



The following statements relate quivers to plabic graphs.

\begin{definition}
Let $G$ be a reduced plabic graph.  We associate a quiver $Q(G)$ as follows.  The vertices of 
$Q(G)$ are labeled by the faces of $G$.  We say that a vertex of $Q(G)$ is \emph{frozen}
if the 
corresponding face is incident to the boundary of the disk, and is \emph{mutable} otherwise.
For each edge $e$ in $G$ which separates two faces, at least one of which is mutable, 
we introduce an arrow connecting the faces;
 this arrow is oriented so that it ``sees the white endpoint of $e$ to the left and the 
black endpoint to the right'' as it crosses over $e$.  We then remove oriented $2$-cycles
from the resulting quiver, one by one, to get $Q(G)$. See 
\cref{fig:plabic2}.
\end{definition}


\begin{lemma}\label{lem:mutG}
If $G$ and $G'$ are related via a square move at a face,
then $Q(G)$ and $Q(G')$ are related via mutation at the corresponding vertex.
\end{lemma}

\acknowledgements{We are grateful to Bernard Leclerc for numerous helpful discussions.
K.S. and M.S.B. acknowledge support from National Science Foundation Postdoctoral Fellowship MSPRF-1502881 and NSF Graduate Research Fellowship
No. DGE-1752814, respectively. 
L. W. was partially supported by 
the NSF grant DMS-1600447.
Any opinions, findings
and conclusions or recommendations expressed in this material are those of 
the authors and do not necessarily reflect the views of the National
Science Foundation. }

\printbibliography

\end{document}